\documentclass[12pt]{article}
\usepackage[cp1251]{inputenc}
\usepackage[english]{babel}
\usepackage{latexsym,amsfonts,amsthm,amsmath,amssymb}

\newtheorem{lemma}{Lemma}[section]
\newtheorem{theorem}[lemma]{Theorem}
\newtheorem{definition}[lemma]{Definition}

\newtheorem{proposition}[lemma]{Proposition}

\newtheorem{example}[lemma]{Example}
\newtheorem{remark}[lemma]{Remark}
\theoremstyle{definition}

\newcommand{\Tor}{\mathop{\rm Tor}\nolimits}

\newcounter{num}

\begin{document}

\title{Complexity and $T$-invariant of Abelian and Milnor groups, and complexity of 3-manifolds}

\author{Ekaterina~{\textsc Pervova}\thanks{Supported by the INTAS YS
fellowship 03-55-1423 and by the Russian Ministry of Education}
\and\addtocounter{footnote}{5} Carlo~{\textsc
Petronio}\footnote{Supported by the INTAS project ``CalcoMet-GT''
03-51-3663}}

\maketitle
\begin{abstract}\noindent We
investigate the notion of complexity for finitely presented groups
and the related notion of complexity for three-dimensional
manifolds.  We give two-sided estimates on the complexity of all
the Milnor groups (the finite groups with free action on $S^3$),
as well as for all finite Abelian groups. The ideas developed in
the process also allow to construct two-sided bounds for the
values of the so-called $T$-invariant (introduced by Delzant) for
the above groups, and to estimate from below the value of
$T$-invariant for an arbitrary finitely presented group. Using the
results of this paper and of previous ones, we then describe an
infinite collection of Seifert three-manifolds for which we can
asymptotically determine the complexity in an exact fashion up to
linear functions.  We also provide similar estimates for the
complexity of several infinite families of Milnor groups.
\vspace{2mm}

\noindent MSC (2000): 20F05 (primary), 57M07, 57M27, 20F65 
(secondary).
\end{abstract}

\renewcommand{\baselinestretch}{1.1}

\section*{Introduction}

The motivation for considering some notion of complexity for
groups is its connection~\cite{MaPe} with the problem of
estimating the complexity of 3-manifolds. The main idea of the
theory of complexity for $3$-manifolds (introduced in~\cite{Ma88,
Ma90}, investigated in~\cite{MaFo}, and comprehensively covered
in~\cite{Ma03}) is to introduce a filtration in the set of all
3-manifolds, in such a way that each level of the filtration
contains only finitely many closed irreducible items.  This allows
to break down the task of classifying all closed 3-manifolds into
an infinite collection of finite classification tasks, because
complexity is also additive under connected sum, so the complexity
of any closed manifold can be computed once the complexity of its
irreducible summands is known. This classification program has
been carried out to a remarkable extent in recent years
(see~\cite{Martelli:surv} and the references quoted therein).

For any given manifold it is very easy to give upper bounds for
its complexity, whereas lower bounds are much harder to establish.
As a matter of fact, the computer programs of Martelli and
Matveev~\cite{Martelli:surv}, which manipulate special spines,
provide upper bounds which experimentally are always sharp. On the
other hand, the only methods currently known to obtain general
lower bounds are those of~\cite{MaPe}, based on group theory (for
some hyperbolic manifolds, there is also a lower estimate in terms
of the volume by Anisov, see \cite{An02}). From this point of view
the results of the present paper can be viewed as potential tools
for constructing more lower estimates on the complexity of
manifolds. This idea is specified in
Lemma~\ref{spine:complexity:lem} below, and a concrete application
is given in Theorem~\ref{compl:mfds:thm}, where we provide
two-sided bounds for the complexity of certain infinite classes of
Seifert manifolds. These estimates are ``asymptotically exact up
to linear maps,'' meaning that the upper and the lower bound
differ by a fixed linear function.

The notion of $T$-invariant, also closely related to the
complexity of 3-manifolds (see~\cite{Del1, Del2, Del3} or the
proof of Theorem~\ref{complexity:Seifert:thm} below), was
introduced by Delzant in~\cite{Del1} for what appear to be
completely different reasons, namely, to study hierarchical
decompositions of finitely presented groups. For instance, it
played a central r\^ole in the proof by Delzant himself and
Potyagailo of the strong accessibility theorem for such
groups~\cite{Del3}. For one-relator groups, the $T$-invariant was
studied in~\cite{KaSch}.

The complexity and the $T$-invariant of a group never coincide,
except for the trivial group, but they are closely related and
they can be studied by similar methods. Exploiting this fact, we
provide in this paper lower bounds for the complexity and the
$T$-invariant of an arbitrary finitely presented group in terms of
the order of the torsion part of its Abelianization
(Theorems~\ref{lower:bound:thm} and~\ref{Delzant:estimate:2:thm}).
Then, for all the members of Milnor's list~\cite{Mil} of finite
groups with free linear action on $S^3$, as well as for all finite
Abelian groups, we present two-sided estimates on their complexity
(Theorems~\ref{upper:Abelian:thm} and~\ref{Milnor:groups:thm}) and
on the value of the $T$-invariant
(Theorems~\ref{Delzant:cyclic:thm} and~\ref{Milnor:groups:thm},
and Remark~\ref{evident:invariant:rem}).

For some Milnor groups, the estimates we obtain show that the
complexity is asymptotically given by the logarithm of the order
of the torsion part of the Abelianization, up to linear maps.
Using results from this paper and from previous ones we also
provide similar ``asymptotically exact'' estimates for the
complexity of certain Seifert manifolds
(Theorem~\ref{compl:mfds:thm}). We then in turn apply this theorem
to derive more precise bounds on the complexity and the
$T$-invariant of some of the other Milnor groups
(Theorems~\ref{Delzant:P:thm} and~\ref{complexity:Seifert:thm}).
The interplay between the estimates on complexity and those on the
$T$-invariant is that by combining upper and lower bounds we can
deduce information on the average length of the relations in a
presentation realizing the complexity
(Propositions~\ref{Petronio:prop},~\ref{Petronio:Milnor:prop},
and~\ref{Petronio:Seifert:prop}). It is interesting to note that
for some Milnor groups the complexity is asymptotically very close
to that of the 3-manifolds whose fundamental groups they are
(Theorems~\ref{compl:mfds:thm} and~\ref{complexity:Seifert:thm}).

\paragraph{Acknowledgments} The authors are grateful
to Sergei~Matveev for his many useful remarks regarding this work.
They would also like to thank Stefan Rosenberg for productive
discussions, Ilya Kapovich and Pierre de la Harpe for introducing
them to Delzant's notion of $T$-invariant, and Nikolai
Moshchevitin for providing important number-theoretical
references. The first-named author would like to thank the
participants of the seminar ``Differential geometry'' directed by
Anatoly Fomenko for useful discussions which led to the refinement
of some ideas used in the paper.

\section{Main definitions}\label{main:def:section}
In this section we define the invariants for which we will provide
estimates in this paper.

\paragraph{Groups}
The notion of group complexity was introduced in~\cite{MaPe}.

\begin{definition}\label{length:def}
{\em Let $\langle a_1, \ldots, a_n|\ r_1, \ldots, r_m \rangle$ be
a presentation of a group. The {\em length} of this presentation
is the number $|r_1|+\ldots+|r_m|$, where $|r_i|$ is the length of
the word $r_i$ in the alphabet $a_1^{\pm 1},\ldots,a_n^{\pm 1}$.
The {\em complexity} $c(G)$ of a group $G$ admitting finite
presentations is the minimum of the lengths of all such
presentations}.
\end{definition}

It can be seen (by explicit enumeration of presentations of small
length) that for $n\leqslant 7$ the complexity of the cyclic group
of order $n$ is equal to $n$. However, the groups $\mathbb Z/_{\!
8}$ and $\mathbb Z/_{\! 9}$ both have complexity $7$, which is
smaller than the order, and $c(\mathbb Z/_{\! {10}})=8$. The
following presentation of $\mathbb Z/_{\! {147}}$, which has
length $23$, shows that the complexity can be significantly
smaller than the order:
$$\langle a,b,c,d|\ a^4bc^4,b^3c^{-1},a^2d^3b^{-1},a^3d^{-1} \rangle.$$

An alternative measure of how complicated a group is, now called
the $T$-invariant of the group, was suggested by Delzant
in~\cite{Del1} and investigated in \cite{Del2,Del3}:

\begin{definition}
{\em The {\em $T$-invariant} $T(G)$ of a finitely presented group
$G$ is the minimal number $t$ such that $G$ admits a presentation
with $t$ relations of length $3$ and an arbitrary number of
relations of length at most $2$. A presentation of this type is
called {\em triangular}.}
\end{definition}

The following easy fact was already noted in~\cite{Del1}:
\begin{proposition}\label{Delzant:def:prop}
$$T(G)=\min\left\{\sum_{i=1}^m\max\{|r_i|-2,0\}:\ G=\langle
a_1,\ldots,a_n|\ r_1,\ldots,r_m \rangle\right\}.$$
\end{proposition}

\paragraph{3-manifolds}
We now review some notions related to 3-manifolds. We will use the
PL category throughout.

\begin{definition}
{\em A 2-dimensional subpolyhedron $P$ of a closed connected
3-manifold $M$ is called a {\em spine} of $M$ if $M\setminus P$ is
homeomorphic to an open 3-ball}.
\end{definition}

In particular, for every spine $P$ of $M$ we have
$\pi_1(P)\cong\pi_1(M)$. We will consider only a particular class
of spines, that we now define.

\begin{definition} \label{simple:def}
{\em A compact polyhedron is called {\em special} if the following
two conditions hold. First, the
link of each point is homeomorphic to one of the following
1-dimensional polyhedra:
\begin{enumerate}
\item[(a)] a circle;
\item[(b)] a circle with a diameter;
\item[(c)] a circle with three radii.
\end{enumerate}
Second, the components of set of points of type (a) are open
discs, while the components of set of points of type (b) are open
segments. The components just described are called {\em faces} and
{\em edges}, respectively, and the points of type (c) are called
{\em vertices}. A {\em special spine} of a closed manifold $M$ is
a spine of $M$ which is a special polyhedron at the same time}.
\end{definition}

The notion of complexity for (arbitrary) 3-manifolds was
introduced in~\cite{Ma88}, see also \cite{Ma90}. We will only need
here the following partial characterization, which could also be
used as a definition:

\begin{proposition}\label{irred:compl:prop}
The {\em complexity} $c(M)$ of a closed irreducible manifold
$M\notin\{S^3,\mathbb R\mathbb P^3,L_{3,1}\}$ is the minimal
number of vertices of a special spine of $M$. The complexity of
the three exceptional manifolds is equal to zero.
\end{proposition}

It turns out that there is a clear relation between the complexity
of a 3-manifold and the complexity of its fundamental group. This
relation is described in the following lemma, which was
essentially proved in~\cite{MaPe}.

\begin{lemma}\label{spine:complexity:lem}
If a manifold $M$ has a special spine $P$ with $n$ vertices then
$\pi_1(M)$ has a presentation of length $3n+3$.
\end{lemma}

\begin{proof}
We know that $\pi_1(M)$ coincides with $\pi_1(P)$. Moreover, the
stratification of $P$ into vertices, edges, and faces gives $P$
the structure of a cell complex. So we can employ the general
algorithm yielding a presentation of the fundamental group of a
cell complex.  The generators are the edges in the complement of a
maximal tree in the $1$-skeleton, so there are $n+1$ of them. The
relations correspond to the faces. Since precisely $3$ faces are
incident to any given edge (with multiplicity), the total length
of the relations is $3(n+1)$.
\end{proof}

According to this result, an upper bound on $c(M)$ implies an
upper bound on $c(\pi_1(M))$, and a lower bound on $c(\pi_1(M))$
implies a lower bound on $c(M)$.

\section{Lower estimates}\label{lower:section}
In this section we establish lower estimates for the complexity
and the $T$-invariant for an arbitrary group, whereas starting
from the next section we will concentrate on Abelian and Milnor
groups.

\paragraph{Group complexity}
We begin with an easy lower bound on the complexity of a group $G$
in terms of the so-called {\em relation matrices} of the
presentations of $G$. Recall that, given a presentation of $G$
with $n$ generators and $m$ relations, the relation matrix
associated to the presentation has size $m\times n$, and its entry
in position $(i,j)$ is the (algebraic) sum of all the exponents of
the $j$-th generator in the $i$-th relation. If $X=(x_{i,j})$ is
the matrix thus obtained we define its {\em norm} as
$$||X||=\sum_{i=1}^m \sum_{j=1}^n|x_{i,j}|.$$
We have the following:

\begin{lemma}\label{evident:lower:lem}
For every finitely presented group $G$ we have
$$c(G)\geqslant \min ||X||,$$
where the minimum is taken over the relation matrices $X$
associated to all possible finite presentations of $G$.
\end{lemma}

\begin{proof}
We only need to note that $||X||$ is less than or equal to the
length of the presentation to which $X$ is associated.
\end{proof}

To provide more effective lower bounds on $c(G)$ we then need to
give an estimate on the possible norms of the relation matrices of
$G$. This is done in the next result, where $\Tor(H)$ denotes the
torsion (\emph{i.e.} finite-order) part of an Abelian group $H$.

\begin{theorem}\label{lower:bound:thm}
For any finitely presented group $G$ we have
$$c(G)\geqslant \log_2|\Tor(G/[G,G])|.$$
\end{theorem}

\begin{proof}
By Lemma~\ref{evident:lower:lem} it is enough to show that
$||X||\geqslant \log_2|\Tor(G/[G,G])|$ for all the relation
matrices $X$ of the presentations of $G$. Fix such an $X$ and
suppose its size is ${m\times n}$. It is well-known that $X$ can
be transformed into a matrix of the form
$$Y=\left(\begin{array}{cc} D & 0 \\ 0 & 0\end{array} \right),$$
with $D$ a $k\times k$ diagonal matrix, and ${\rm det}(D)\ne 0$ if
$k\geqslant 1$, using a finite sequence of operations as follows:
\begin{enumerate}
\item Interchange two rows or columns; \item Multiply one row or
column by $-1$; \item Add one row or column to a different one.
\end{enumerate}
One can now see that each such operation transforms the relation
matrix of a finite presentation of $G$ into the relation matrix of
some other presentation of $G$, which easily implies that
$|\Tor(G/[G,G])|=|{\rm det}(D)|$.

If we set $d=|{\rm det}(D)|$ we have the obvious property that
{\em the determinants of all the $k\times k$ submatrices of $Y$
are multiples of $d$, and some of them is non-zero}. Moreover one
can easily see that this property is preserved under all the
inverse operations which lead from $Y$ back to $X$. Therefore
there is a $k\times k$ submatrix $X'$ of $X$ such that ${\rm
det}(X')$ is a non-zero multiple of $d$, so in particular $|{\rm
det}(X')|\geqslant d$.

We can now observe that $|{\rm det}(X')|$ is bounded from above by
the product of the Euclidean norms of the rows of $X'$, and each
such norm is bounded from above by the $L^1$-norm, whence by the
$L^1$-norm of the corresponding whole row of $X$. Noting that each
non-zero row of $X$ has norm at least $1$, and dismissing the zero
rows if necessary (recall that we want to give a lower bound on
$\|X\|$), we conclude that $|{\rm det}(X')|$ is bounded from above
by the product of the $L^1$-norms of all the non-zero rows of
$X$. Therefore
\begin{eqnarray*}
&&|\Tor(G/[G,G])|=d\leqslant \prod_{i=1}^m\sum_{j=1}^n|x_{i,j}|\\
&\Rightarrow&
\log_2|\Tor(G/[G,G])|\leqslant\sum_{i=1}^m\log_2\left(\sum_{j=1}^n|x_{i,j}|\right).
\end{eqnarray*}
Noting that $\log_2n\leqslant n$ for all $n\in\mathbb N$ we deduce
that
$$\log_2|\Tor(G/[G,G])|\leqslant\sum_{i=1}^m\sum_{j=1}^n|x_{i,j}|=\|X\|,$$
and the proof is complete.
\end{proof}

\begin{remark}\label{first:improve:base:rem}
{\em In the previous statement the base $a=2$ of logarithms was
chosen because it has the property that $\log_a(n)\leqslant n$ for
all $n\in\mathbb N$, and the theorem remains true with any other
base $a$ having this property. One easily sees that the best lower
estimate for $c(G)$ is obtained for $a=\sqrt[3]3$. Since we are
only interested in the qualitative fact that a logarithmic lower
bound exists, we will keep employing the base $2$. However we will
use the fact that the inequality in the previous statement is
strict}.
\end{remark}

\begin{remark}\label{real:lower:bound:rem}
{\em Along the proof of Theorem~\ref{lower:bound:thm} we have
shown that for any presentation $\langle a_1, \ldots, a_n|\ r_1,
\ldots, r_m \rangle$ of a group $G$ the following inequality is
valid:}
$$\log_2|\Tor(G/[G,G])|\leqslant\log_2|r_1|+\ldots+\log_2|r_m|.$$
\end{remark}

\paragraph{The $T$-invariant}
Proposition~\ref{Delzant:def:prop} allows to conclude immediately
that $T(G)< c(G)$ for a non-trivial $G$. However, we will show
that in many cases the invariants $c$ and $T$ are asymptotically
equivalent. We begin with two rather easy results.

\begin{lemma}\label{Delzant:pres:lem}
Let $G$ be a finitely presented group. Then for every finite
presentation of $G$ there is a presentation of the same or smaller
length which contains relations of length $\geqslant 2$ only, and
all relations of length $2$ are of the form $x^2$. Moreover, $G$
admits a triangular presentation with exactly $T(G)$ relations of
length $3$ and some relations of the form $x^2$.
\end{lemma}

\begin{proof}
Of course the relations of length $0$ can always be omitted.
Suppose the given presentation is $\langle a_1,\ldots, a_n|\
r_1,\ldots,r_t,r_1',\ldots,r_k' \rangle$, where the lengths of all
$r_i$ are $\geqslant 3$, and the lengths of all $r_j'$ are $1$ or
$2$. We describe a recursive procedure to get a presentation as
desired. If some $r_j'$ has length $1$, {\em i.e.} it is of the
form $x$, then this $x$ can be removed from the generators and
from all the relations where it occurs. This produces a
presentation of $G$ with the same or a smaller number of relations
of length $\geqslant 3$, and the total length of the presentation,
as well as each relation's own length, is not increased. Hence
there is a presentation of $G$ of the same or smaller length with
$\leqslant t$ relations of length $\geqslant 3$ and some number of
relations of length $2$.

Consider a relation of length $2$. If it has the form $x^{-1}y$,
then we can discard it and the generator $y$, replacing all
occurrences of $y$ in all the relations by $x$. Such a procedure
increases neither the lengths of the other relations nor the total
length of the presentation. So in the end we get a presentation of
$G$ of the described type and of the same or smaller length.

If the initial presentation is triangular, so is the final one,
and the last assertion follows.
\end{proof}

\begin{proposition}\label{Delzant:estimate:1:prop}
Let $G$ be a nontrivial finitely presented group without
$2$-torsion. Then
$$\frac{1}{3}c(G)\leqslant T(G)< c(G).$$
\end{proposition}

\begin{proof}
The second inequality is valid in general and was already remarked
above. For the first inequality, consider a triangular
presentation of $G$ of the type described in
Lemma~\ref{Delzant:pres:lem}. Suppose there is a relation of the
form $x^2$. Then, since $G$ does not have elements of order $2$,
$x$ has to be trivial in $G$. Hence we can remove all occurrences
of $x$ from all the relations and remove $x$ itself from the list
of generators. As a result we get a triangular presentation of $G$
which contains exactly $T(G)$ relations of length $3$ (a priori no
more than that number, but by minimality there cannot be less) and
no relations of length less than $3$.
\end{proof}

Combining the previous result with Theorem~\ref{lower:bound:thm}
we deduce that if $G$ has no $2$-torsion then
$$T(G)\geqslant\frac 13 \log_2|\Tor(G/[G,G])|.$$
This estimate can actually be improved to a stronger and general
one:

\begin{theorem}\label{Delzant:estimate:2:thm}
Let $G$ be a finitely presented group, and let
$|\Tor(G/[G,G])|=2^l(2m+1)$. Then
$$T(G)\geqslant \log_3(2m+1).$$
\end{theorem}

\begin{proof}
Consider a triangular presentation of the type described in
Lemma~\ref{Delzant:pres:lem}, with $t=T(G)$ relations of length
$3$ and some $h$ relations of the form $x^2$. The relation matrix
$X$ of this presentation has the form
$$X=\left(\begin{array}{cc} Y & Z \\ 2I_h & 0\end{array}\right),$$
where the matrix $(Y\, Z)$ has $t$ rows, the $L^1$-norm of each
row is $3$, and $I_h$ is the $h\times h$ unit matrix. As in the
proof of Theorem~\ref{lower:bound:thm}, $X$ has a square submatrix
$X'$ whose determinant is non-zero and divisible by
$|\Tor(G/[G,G])|=2^l(2m+1)$. On the other hand, $X'$ consists of
$t'\leqslant t$ rows of $L^1$-norm $3$ and some $h'\leqslant h$
rows with a single non-zero entry equal to $2$. It follows that
$|{\rm det}(X')|=2^{h'}\cdot\delta$ and $\delta\leqslant
{3}^{t'}\leqslant 3^{t}$. Recalling that $|{\rm det}(X')|$ is
divisible by $2^l(2m+1)$, we deduce that $\delta$ is divisible by
$2m+1$, whence $2m+1\leqslant\delta\leqslant 3^t$, and the desired
estimate follows.
\end{proof}

\section{Abelian and Milnor groups}\label{upper:geom:section}
In this section we will give two-sided estimates for the
complexity and the $T$-invariant of all the Milnor groups, and in
the process we will obtain similar estimates for Abelian groups.

\paragraph{Milnor groups}
As already mentioned, Milnor classified in~\cite{Mil} the finite
groups having a free linear action on $S^3$. The complete list of
all such groups is as follows:
\begin{enumerate}
\item $Q_{4n}=\langle x,y|\ x^{-1}yxy,\, x^{-2}y^n\, \rangle$,
$n\geqslant 2$; \item $D_{2^k(2n+1)}=\langle x,y|\ x^{2^k},\,
y^{2n+1},\, xyx^{-1}y\,\rangle$, $k\geqslant 3$, $n\geqslant 1$;
 \item
 \[\left.
 \begin{array}{l}
 P_{24}=\langle x,y|\ x^{-1}yxyxy,\, x^{-2}y^3,\, x^4\, \rangle, \\
 P_{48}=\langle x,y|\ x^{-1}yxyxy,\, x^{-2}y^4,\, x^4\, \rangle, \\
 P_{120}=\langle x,y|\ x^{-1}yxyxy,\, x^{-2}y^5,\, x^4\, \rangle,
 \\
 P_{8\cdot 3^k}'=\langle x,y,z|\ x^{-1}yxy,\, x^{-2}y^2,\, zxz^{-1}y^{-1},\, zyz^{-1}y^{-1}x^{-1},\,
 z^{3^k}\, \rangle,\;\; k\geqslant 2; \\
 \end{array}
 \right.\]
 \item The direct product of any of the groups listed so far, or of the trivial group, with a cyclic group
 of coprime order.
\end{enumerate}

We will start by considering the simplest case, namely, cyclic
groups. It is interesting to note that many of the ideas which
will be used for the other groups in the list are already present
at this level.

\paragraph{Cyclic groups}
In this paragraph we give two-sided estimates for the complexity
and the $T$-invariant of finite cyclic groups and we describe some
properties of their minimal presentations. Although arbitrary
finite Abelian groups are not in Milnor's list, the complexity
estimates can be generalized to include them with no extra effort,
so we cover them too. We begin with an upper estimate:

\begin{proposition}\label{upper:cyclic:prop}
For every $p\geqslant 2$ we have $c(\mathbb Z/_{\! p})< 4\log_2p$.
\end{proposition}

To prove this result, we first state the following evident:

\begin{lemma}\label{division:lemma}
Let $p,q,r,s$ be non-negative integers such that $p=sq+r$.
Consider a group presentation of length $l$ involving a relation
of the form $ua^pv$, where $a$ is a generator and $u,v$ are words
in the generators. Then a new presentation of the same group is
obtained by adding a new generator $b$ and a new relation
$b^{-1}a^q$, and replacing the relation $ua^pv$ by $ub^sa^rv$. The
length of the new presentation is $1+l+q+r+s-p$.
\end{lemma}

Let us denote now by $\ell(p)$ the shortest length of a
presentation of the group $\mathbb Z/_{\! p}$ obtained from the
trivial presentation $\langle a|\ a^p\rangle$ by repeated
application of the lemma just stated. Of course $c(\mathbb Z/_{\!
p})\leqslant \ell(p)$, so the next result implies
Proposition~\ref{upper:cyclic:prop}:

\begin{proposition}\label{ell:prop} For all $p\geqslant 2$ we have
$\ell(p)< 4 \log_2p$.
\end{proposition}

\begin{proof}
We proceed by induction on $p$, noting that the inequality is true
for $p=2$ and $p=3$, because $\ell(2)=2$ and $l(3)=3$. For the
inductive step, for $p\geqslant 4$ we apply
Lemma~\ref{division:lemma} with $q=2$, and we distinguish
according to the parity of $p$. If $p$ is even then
\begin{eqnarray*}
&\mathbb Z/_{\!
p}&=\langle a,b|\ b^{-1}a^2,b^{p/2}\rangle\\
&\Rightarrow& \ell(p)\leqslant
3+\ell(p/2)<3+4\log_2(p/2)=4\log_2p-1<4\log_2(p).
\end{eqnarray*} If $p$ is odd then
\begin{eqnarray*}
&\mathbb Z/_{\!
p}&=\langle a,b|\ b^{-1}a^2,b^{(p-1)/2}a\rangle\\
&\Rightarrow& \ell(p)\leqslant 4+\ell((p-1)/2)\\
& & \qquad <4+4\log_2((p-1)/2)=4\log_2(p-1)<4\log_2(p).
\end{eqnarray*} This proves the
desired inequality.
\end{proof}

\begin{example}
{\em Applying the procedure described in the previous proof for
$p={357}$ we get the following length-27 presentation of $\mathbb
Z/_{\! {357}}$:}
$$\langle a,b,c,d,e,f,g,h|\
b^{-1}a^2,c^{-1}b^2,d^{-1}c^2,e^{-1}d^2,f^{-1}e^2,g^{-1}f^2,h^{-1}g^2,h^2gfca
\rangle.$$
\end{example}

\begin{remark}\label{second:improve:base:rem}
{\em The estimate in Proposition~\ref{ell:prop}, whence that in
Proposition~\ref{upper:cyclic:prop}, can actually be improved
using $q=3$ rather than $q=2$ for the repeated application of
Lemma~\ref{division:lemma}. Namely, one can show that
$c(\mathbb{Z}/_{\! p})< 6 \log_3p,$ which is a slightly better
bound since $\frac{6}{\log_23}\approx 3.786< 4$. All other choices
of $q$, on the other hand, produce bounds with larger constants.}
\end{remark}

In the case of an arbitrary finite Abelian group $G$,
Proposition~\ref{upper:cyclic:prop} implies the following result.
Recall that the rank of $G$ is the minimal number of cyclic groups
which $G$ can be expressed as the direct product of.

\begin{theorem}\label{upper:Abelian:thm}
Let $G$ be a finite Abelian group of rank $k$. Then $c(G)<
4\log_2|G|+2k(k-1)$.
\end{theorem}

\begin{proof}
By assumption, $G$ is isomorphic to $\mathbb Z/_{\!
{p_1}}\oplus\ldots\oplus \mathbb Z/_{\! {p_k}}$, and
$|G|=p_1\cdot\ldots \cdot p_k$. We can obtain a presentation of
$G$ by taking the union of the presentations of $\mathbb Z/_{\!
{p_i}}$ constructed in Proposition~\ref{ell:prop} and adding
relations providing the commutativity. Notice that the generating
set constructed in Proposition~\ref{ell:prop} always contains an
element generating the whole group. Hence it suffices to add
$k(k-1)/2$ relations of length $4$ each (the commutators of all
pairs of different generators). This produces a presentation of
length that is strictly less than
\begin{eqnarray*}
& & 4\log_2p_1+\ldots+4\log_2p_k+4\frac{k(k-1)}{2}\\
& = & 4\log_2(p_1\cdot\ldots \cdot p_k)+2k(k-1)=4\log_2|G|+2k(k-1)
\end{eqnarray*}
and the theorem is proved.
\end{proof}

\paragraph{Average length}

In view of Proposition~\ref{Delzant:def:prop}, in order to
estimate the $T$-invariant for a cyclic group, one would need a
lower bound on the number of relations in a length-minimizing
presentation ({\em i.e.} a presentation realizing the complexity).
Such a bound is established in the following:

\begin{proposition}\label{Petronio:prop}
Suppose that $\langle a_1,\ldots,a_n|\ r_1,\ldots,r_m \rangle$ is
a length-minimizing presentation of a finite cyclic group. Then
$$\frac{|r_1|+\ldots+|r_m|}{m}< 16.$$
\end{proposition}

\begin{proof}
If $p$ is the order of $G$, according to
Proposition~\ref{upper:cyclic:prop} we have
$\log_2p>(|r_1|+\ldots+|r_m|)/4$. Moreover $p= |\Tor(G/[G,G])|$,
so by Remark~\ref{real:lower:bound:rem} we have
$\log_2p\leqslant\log_2|r_1|+\ldots+\log_2|r_m|$. Combining these
two estimates we get the inequality
$$|r_1|+\ldots+|r_m|<4\cdot\big(\log_2|r_1|+\ldots+\log_2|r_m|\big),$$
which we divide by $m$ to obtain
$$\frac{|r_1|+\ldots+|r_m|}{m}<4\frac{\log_2|r_1|+\ldots+\log_2|r_m|}{m}.$$
The right-hand side of the latter inequality is just
$4\log_2\sqrt[m]{|r_1|\cdot\ldots\cdot|r_m|}$. Applying the Cauchy
inequality between the geometric mean and the arithmetic one to
the expression under the sign of logarithm, we get
$$\frac{|r_1|+\ldots+|r_m|}{m}<4\log_2\left(\frac{|r_1|+\ldots+|r_m|}{m}\right).$$
It follows that the number $(|r_1|+\ldots+|r_m|)/m$ must satisfy
the inequality $x< 4\log_2x$. Now, the numbers satisfying this
inequality lie between the two solutions to the equation
$x=4\log_2x$ and therefore are bounded from above by the greatest
of them, which evidently is $x=16$, because the other solution
lies between $1$ and $2$. The desired estimate on
$(|r_1|+\ldots+|r_m|)/m$ follows.
\end{proof}

\begin{remark}
{\em Using the better estimate given by
Remark~\ref{second:improve:base:rem} one could show that
$(|r_1|+\ldots+|r_m|)/m<15$.}
\end{remark}

The values of the $T$-invariant for cyclic groups can now be
estimated as follows.

\begin{theorem}\label{Delzant:cyclic:thm}
For every odd $p$ and for every integer $k\geqslant 0$ we have
that
$$\frac{1}{\log_23}\log_2p\;\leqslant\; T\left(\mathbb{Z}/_{\! {2^kp}}\right)\;<\; \frac{7}{2}(\log_2p+k).$$
\end{theorem}

\begin{proof}
The lower estimate follows from
Theorem~\ref{Delzant:estimate:2:thm}. For the upper estimate, take
a length-minimizing presentation of $\mathbb{Z}/_{\! {2^kp}}$ with
relations $r_1,\ldots,r_m$. Since $|r_i|\geqslant 2$,
Proposition~\ref{Delzant:def:prop} implies that
$$T(\mathbb{Z}/_{\! {2^kp}})\leqslant |r_1|+\ldots+|r_m|-2m.$$

Proposition~\ref{Petronio:prop} now yields
$$|r_1|+\ldots+|r_m|-2m\leqslant \frac 78 (|r_1|+\ldots+|r_m|),$$
and $|r_1|+\ldots+|r_m|=c\left(\mathbb{Z}/_{\! {2^kp}}\right)$ by
the choice of the presentation, but $c\left(\mathbb{Z}/_{\!
{2^kp}}\right)<4(\log_2p+k)$ by
Proposition~\ref{upper:cyclic:prop}, whence the conclusion.
\end{proof}

\paragraph{Other Milnor groups}
A straight-forward application of the technique used to prove
Proposition~\ref{ell:prop} and of Theorem~\ref{lower:bound:thm}
allows us to obtain some estimates on the complexity of all the
other groups in Milnor's list:

\begin{theorem}\label{Milnor:groups:thm}
The following estimates hold for the
complexity of the Milnor groups, where in all cases the lower
bound equals the base-$2$ logarithm of the order of the torsion of
the Abelianization, and the same term appears in the upper
estimate too:
\begin{enumerate}
\item For every $n\geqslant 2$ and every odd $q$ coprime with $n$
we have
$$\log_2q+2<c(Q_{4n}\times \mathbb{Z}/_{\! q})<
4(\log_2q+2)+4\log_2n+6;$$
 \item For every $k\geqslant 3$ and every coprime odds $n\geqslant 3,q\geqslant 1$,
 we have
 $$\log_2q+k <c(D_{2^kn}\times\mathbb{Z}/_{\! q}) < 4(\log_2q+k)+4\log_2n+12;$$
 \item For every $q$ coprime with $2$ and $3$ (and $5$, for the last estimate) we have
 \[\left.
 \begin{array}{l}
 \log_2(3q)< c(P_{24}\times\mathbb{Z}/_{\! q})< 4(\log_2(3q))+17; \\
 \log_2q+1<c(P_{48}\times\mathbb{Z}/_{\! q})< 4(\log_2q+1)+20; \\
 \log_2q< c(P_{120}\times\mathbb{Z}/_{\! q})< 4\log_2q+25; \\
 \end{array}
 \right.\]
 \item For every $k\geqslant 2$ and every $q$ coprime with $2$ and $3$
 we have
 $$\log_2q+(\log_23)k<c(P_{8\cdot 3^k}'\times\mathbb{Z}/_{\! q})< 4(\log_2q+(\log_23)k)+29.$$
\end{enumerate}
\end{theorem}

\begin{proof} The lower bounds are obtained by direct application
of Theorem~\ref{lower:bound:thm} and
Remark~\ref{first:improve:base:rem}. To get the upper bounds, we
start by writing the most straight-forward presentation for each
of the groups listed in the theorem. Namely, we add to each of the
standard presentations reproduced at the beginning of the present
section, one generator $a$ (corresponding to $\mathbb{Z}/_{\!
q}$), a relation $a^q$, and the commutation relations $[x,a]$ and
$[y,a]$ (and $[z,a]$ for the group $P_{8\cdot 3^k}'$).

Now we apply exactly the trick described in
Lemma~\ref{division:lemma}. Evidently, in case of $Q_{4n}$ this
produces a presentation of length $\ell(n)+\ell(q)+14$. From the
estimate given by Proposition~\ref{ell:prop}, we deduce that this
number is less than $4\log_2(nq)+14$, and the conclusion easily
follows. In case of $D_{2^kn}$ the trick produces a presentation
of length $\ell(2^k)+\ell(n)+\ell(q)+12$, which implies the bound
stated in the theorem.  For each of $P_{24}$, $P_{48}$, $P_{120}$
we can get a presentation of that group multiplied by
$\mathbb{Z}/_{\! q}$ having length $\ell(q)$ plus the length of
the presentation of that group given above plus $8$, and the upper
bounds follow after easy calculations. Finally, for $P_{8\cdot
3^k}'$ we get a presentation of length $\ell(3^k)+\ell(q)+29$,
which gives the desired bound again.
\end{proof}

\begin{remark}
{\em Slightly better numerical estimates could be shown using
Remarks~\ref{first:improve:base:rem}
and~\ref{second:improve:base:rem}.}
\end{remark}

\begin{remark}\label{evident:invariant:rem}
{\em The above result together with
Theorem~\ref{Delzant:estimate:2:thm} and the inequality
$T(G)<c(G)$ allows to obtain upper and lower bounds on the
$T$-invariant of all Milnor groups. Since the upper bounds thus
obtained coincide with those for complexity, and the lower bounds
are just an immediate consequence of
Theorem~\ref{Delzant:estimate:2:thm}, we do not spell them out
here (see also below).}
\end{remark}

The complexity estimates given in
Theorem~\ref{Milnor:groups:thm}(3,4) are ``asymptotically exact up
to linear functions,'' and we can now exploit this fact to give a
similar estimate also for the $T$-invariant.  We begin with the
following:

\begin{proposition}\label{Petronio:Milnor:prop}
Let $P$ be one of the groups $P_{24}$, $P_{48}$, $P_{120}$,
$P_{8\cdot 3^k}'$, and let $q$ be a positive integer coprime with
$|P|$. Suppose that $\langle a_1,\ldots, a_n|\ r_1,\ldots, r_m
\rangle$ is a length-minimizing presentation of
${P}\times\mathbb{Z}/_{\! q}$. Then
$$\frac{|r_1|+\ldots+|r_m|}{m}< 52.$$
\end{proposition}

\begin{proof}
Let $G$ denote ${P}\times\mathbb{Z}/_{\! q}$. Combining the lower
estimate given by Remark~\ref{real:lower:bound:rem} with the upper
estimate given by Theorem~\ref{Milnor:groups:thm}(3,4), we deduce
that
$$|r_1|+\ldots+|r_m|<4\log_2\big(|r_1|\cdot\ldots\cdot|r_m|\big)+29.$$
Dividing by $m$, noting that $29/m\leqslant 29$, and using the
Cauchy inequality as in the proof of
Proposition~\ref{Petronio:prop}, we deduce that
$$\frac{|r_1|+\ldots+|r_m|}{m}< 4\log_2\left(\frac{|r_1|+\ldots+|r_m|}{m}\right)+29,$$
which means that the number $(|r_1|+\ldots+|r_m|)/m$ must satisfy
the inequality $x< 4\log_2x+29$. Any such $x$ lies between the two
solutions of the equation $x=4\log_2x+29$. Notice that the smaller
of the two solutions lies between $\frac{1}{2^8}$ and
$\frac{1}{2^7}$, and that $52$ does not satisfy the inequality $x<
4\log_2x+29$. It follows that all solutions to this inequality are
less than $52$.
\end{proof}

Once again, we could slightly improve the numerical estimate given
by this proposition, but we are actually only interested in the
fact that a fixed upper bound exists. As announced, we use the
proposition to give asymptotically exact estimates for
$T$-invariant.

\begin{theorem}\label{Delzant:P:thm}
Let ${P}$ be one of the groups $P_{24}$, $P_{48}$, $P_{120}$. Let
$q$ be a positive integer coprime with $2$ and $3$ (and $5$, for
$P=P_{120}$).  Then
$$\frac{1}{\log_23}\log_2q\leqslant T({P}\times\mathbb{Z}/_{\! q})< \frac{50}{13}\log_2q+24,$$
$$\frac{1}{\log_23}\log_2q+k \leqslant T(P_{8\cdot
3^k}'\times\mathbb{Z}/_{\! q})< \frac{50}{13}(\log_2q+\log_23\cdot
k)+29.$$
\end{theorem}

\begin{proof}
Let $G$ be one of the groups mentioned in the statement. Since $q$
is odd, the lower bounds follow from
Theorem~\ref{Delzant:estimate:2:thm}. To get the upper bounds, we
again apply Proposition~\ref{Delzant:def:prop} to a
length-minimizing presentation of $G$, getting
$$T(G)\leqslant c(G)-2m,$$
where $m$ is the number of relations of a length-minimizing
presentation of $G$. Now Proposition~\ref{Petronio:Milnor:prop}
implies that $m>\frac{c(G)}{52}$, so we finally get
$$T(G)<\frac{25}{26}c(G).$$
Combining this inequality with the upper bounds from
Theorem~\ref{Milnor:groups:thm}, we get our statement.
\end{proof}

\section{More asymptotically exact estimates}
For the Milnor groups of type $P_*^{(')}\times\mathbb{Z}/_{\! q}$,
Theorem~\ref{Milnor:groups:thm}(3,4) provides estimates which are
``asymptotically exact up to linear functions'' as $q$ (and $k$)
tend to infinity, because the upper and lower estimates only
differ by a fixed linear function. This is not the case for the
other Milnor groups, {\em i.e.} those of
Theorem~\ref{Milnor:groups:thm}(1,2), but it turns out that there
are infinite families of such groups for which similar asymptotic
estimates actually do hold. This section is devoted to these
estimates and to some related ones, having the same property of
``asymptotic exactness,'' for the $T$-invariant of the same groups
and for the complexity of certain Seifert $3$-manifolds. As a
matter of fact, the estimates for groups depend on those for
$3$-manifolds, which employ results established elsewhere by more
geometric methods.

\paragraph{Zaremba pairs}
To describe the families of Milnor groups we will deal with we
must make a digression into number theory. Specifically, we need
the following definition and some facts related to it.

\begin{definition}\label{Zaremba:pair:def}
{\em A pair of coprime positive integer numbers $(p,q)$ with
$p>q$, is called a {\em Zaremba pair} if all the partial quotients
$a_i$ in the expansion of $p/q$ into the continuous fraction
$$\frac{p}{q}=a_1+\frac{1}{a_2+\frac{1}{\ldots+\frac{1}{a_n}}}$$
satisfy the inequality $a_i\leqslant 5$}.
\end{definition}

For example, all pairs of consecutive Fibonacci numbers are
Zaremba pairs.

For all coprime $p>q>1$ we denote now by $S(p,q)$ the sum of all
the partial quotients in the expression of $p/q$ as a continued
fraction.

\begin{proposition}\label{desired:bound:prop}
If $(p,q)$ is a Zaremba pair then $S(p,q)\leqslant 3\log_2p$.
\end{proposition}

\begin{proof}
For any coprime $p'>q'>1$, if $a_1,\ldots,a_n$ are the partial
quotients in the expression of $p'/q'$ as a continued fraction,
one can easily show by induction that
$$\left(
\begin{array}{c}
p' \\
q'
\end{array}
\right) = \left(
\begin{array}{cc}
a_1 & 1 \\
1 & 0
\end{array}
\right) \left(
\begin{array}{cc}
a_2 & 1 \\
1 & 0
\end{array}
\right) \cdots \left(
\begin{array}{cc}
a_{n-1} & 1 \\
1 & 0
\end{array}
\right) \left(
\begin{array}{c}
a_n \\
1
\end{array}
\right). $$ For a Zaremba pair $(p,q)$ only five matrices can
appear in this formula
$$\label{matrices:allowed:formula} \left(
\begin{array}{cc}
1 & 1 \\
1 & 0
\end{array}
\right),\quad \left(
\begin{array}{cc}
2 & 1 \\
1 & 0
\end{array}
\right),\quad \left(
\begin{array}{cc}
3 & 1 \\
1 & 0
\end{array}
\right),\quad \left(
\begin{array}{cc}
4 & 1 \\
1 & 0
\end{array}
\right),\quad \left(
\begin{array}{cc}
5 & 1 \\
1 & 0
\end{array}
\right), $$
and there are only four possible ``starting points''
$\left(
\begin{array}{c}
a_n \\
1
\end{array}
\right),$ namely
$$\left(
\begin{array}{c}
2 \\
1
\end{array}
\right),\qquad \left(
\begin{array}{c}
3 \\
1
\end{array}
\right),\qquad \left(
\begin{array}{c}
4 \\
1
\end{array}
\right),\qquad \left(
\begin{array}{c}
5 \\
1
\end{array}
\right).$$

The proof now proceeds by induction on the length $n$ of the
expansion. For $n=1$ the conclusion follows from the fact that
that $m\leqslant 3\log_2m$ for $2\leqslant m\leqslant 5$.

For the inductive step we note that if $2\leqslant m\leqslant 5$ and we
are given $(p',q')$ such that $S(p',q')\leqslant 3\log_2p'$,
setting
$$\left(\begin{array}{c} p'' \\ q''
\end{array}\right)=\left(\begin{array}{cc} m & 1 \\ 1 &
0\end{array}\right) \left(\begin{array}{c} p' \\ q'
\end{array}\right),$$
we have
\begin{eqnarray*}
S(p'',q'') & = & m+S(p',q')\leqslant
3\log_2m+3\log_2p'\\
& = & 3\log_2(mp')\leqslant 3\log_2(mp'+q')=3\log_2p''.
\end{eqnarray*}
This does not quite suffice to conclude when the expansion
involves some matrix with 1 in position (1,1). However one
notes that the inequality $m\leqslant 3\log_2m$ actually holds also for $m=6$.
Therefore, if one of $m',m''$ is 1 and the other one is between 1 and 5,
given $(p',q')$ such that
$S(p',q')\leqslant 3\log_2p'$, if we set
$$\left(\begin{array}{c} p'' \\ q''
\end{array}\right)=\left(\begin{array}{cc} m' & 1 \\ 1 &
0\end{array}\right) \left(\begin{array}{cc} m'' & 1 \\ 1 &
0\end{array}\right) \left(\begin{array}{c} p' \\ q'
\end{array}\right)$$
we have
\begin{eqnarray*}
S(p'',q'') & = & m'+m''+S(p',q')\leqslant
3\log_2(m'+m'')+3\log_2p'\\
& = & 3\log_2((m'+m'')p')\leqslant 3\log_2((m'+m'')p'+m'q')=3\log_2p''.
\end{eqnarray*}
This argument suffices to prove the inequality for all Zaremba pairs
except those of type $(a+1,a)$, for which the conclusion is obvious.
\end{proof}

The conditions of Definition~\ref{Zaremba:pair:def} may appear to
be rather restrictive. However, the following fact is conjectured
by numerical analysts: {\em there exists a constant $B$ with the
property that for every $p$ there exists $1<q<p$ coprime with $p$
such that all partial quotients in the expansion of $p/q$ as a
continued fraction are not greater than $B$}. (This statement is
known as {\em Zaremba's conjecture}, see~\cite{Za72}. Its
motivation is to find optimal lattice points for numerical
integration, see also \cite{Kor}). Cusick conjectured
in~\cite{Cu93} that $B=5$. So far Zaremba's conjecture has been
proved only in a few particular cases. Niederreiter proved it for
powers of $2$ and $3$~\cite{Nie86}, and Yodphotong and Laohakosol
proved it for powers of $6$~\cite{LaYo}. On the other hand, it is
known that there are actually ``many'' Zaremba
pairs~\cite{He2,He1}. We will also use the following weaker
definition.

\begin{definition}
{\em A pair of coprime positive integers $(p,q)$ with $p>q$ is
called a {\em weak Zaremba pair} if the partial quotients
$a_1,\ldots,a_n$ in the expansion of $p/q$ into continuous
fraction satisfy the inequality ${a_1+\ldots+a_n}\leqslant 5n$}.
\end{definition}

Weak Zaremba pairs were investigated in \cite{Cu85,He3,Co03},
where it was shown that they are also not infrequent.

\begin{proposition}\label{weak:bound:prop}
If $(p,q)$ is a weak Zaremba pair then $S(p,q)\leqslant
10\log_2p$.
\end{proposition}

\begin{proof}
An easy induction argument shows that $n\geqslant 2\log_2p$,
whence the conclusion at once.
\end{proof}

\paragraph{Asymptotically exact estimates for manifolds}
In this paragraph we consider the complexity of certain Seifert
manifolds and of their fundamental groups, which in some cases
allows to obtain better bounds than those provided by
Theorem~\ref{Milnor:groups:thm}. We employ for Seifert manifolds
the same notation as in~\cite{MarPet}. Namely, if $F$ is a closed
surface, $t$ is an integer, and $(p_1,q_1),\ldots,(p_k,q_k)$ are
coprime pairs of integers with $|p_i|\geqslant 2$, then
$$\big(F;(p_1,q_1), \ldots,(p_k,q_k),t\big)$$ denotes the (oriented) Seifert
manifold obtained from $F\times S^1$ or from $F\tilde{\times}
S^1$, if $F$ is nonorientable, by removing $k+1$ solid fibred tori
and performing Dehn filling on the resulting boundary components
with slopes $p_1a_1+q_1b_1,
\ldots,p_ka_k+q_kb_k,a_{k+1}+tb_{k+1}$. Here the $a_i$'s are
contained in a section of the bundle, the $b_i$'s are fibres, and
each $a_i,b_i$ is a positive basis in homology.

\begin{theorem}\label{compl:mfds:thm}
Let $(p,q)$ be a Zaremba pair. Then
$$\frac{2}{\log_25}\log_2p-1\leqslant c(L_{p,q})\leqslant 3\log_2p-3,$$
$$\frac{2}{\log_25}\log_2q\leqslant c(S^2;(2,1),(2,1),(p,q),-1) < 3\log_2q+9.$$
\end{theorem}

\begin{proof}
The proof of the upper bound in the first formula is by direct
application of~\cite[Theorem~1.11]{MarPet} and
Proposition~\ref{desired:bound:prop} above, because, with the
notation of~\cite{MarPet}, we have $S(p,q)=|p,q|+1$. For
the manifold $M$ of the second
formula,~\cite[Theorem~1.11]{MarPet}
and~\cite[Theorem~2.5]{MarPet} yield the bound
$$c(M)\leqslant S(p,q)+1.$$
Since $(p,q)$ is a Zaremba pair,
Proposition~\ref{desired:bound:prop} implies that $S(p,q)\leqslant
3\log_2p$. However, as for any Zaremba pair, $p/q<a_1+1\leqslant
6$. Hence
$$c(M)< 3\log_2q+3\log_26+1=3\log_2q+\log_227+4.$$
Since $c(M)$ is an integer, we conclude that it actually does not
exceed $3\log_2q+([\log_227]+1)+4=3\log_2q+9$, as desired.

The lower bounds in both formulae follow directly
from~\cite[Theorem~1]{MaPe}, because the order of the first
homology group is $p$ for $L_{p,q}$, and it is $4q$ for the
manifold in the second formula. This proves the theorem.
\end{proof}

\begin{remark}
{\em If $(p,q)$ is a weak Zaremba pair and $M$ is the lens space
with parameters $(p,q)$ then the complexity of $M$ still depends
on the order of the first homology group of $M$ logarithmically,
as in point 1 of the previous theorem, except that the constants
are worse. Namely, we have}
$$\frac{2}{\log_25}\log_2p-1\leqslant c(M)\leqslant 10\log_2p-3.$$
\end{remark}

\paragraph{Exact asymptotic estimates for groups}
Theorem~\ref{compl:mfds:thm} serves as a tool to obtain good
estimates on the complexity and the $T$-invariant for some of the
Milnor groups in Theorem~\ref{Milnor:groups:thm}(1,2).

\begin{theorem}\label{complexity:Seifert:thm}
\begin{enumerate}
\item Let $(n,q)$ be a Zaremba pair with odd $q$. Then
$$\begin{array}{rcccl}
\log_2q+2 & < & c(Q_{4n}\times \mathbb{Z}/_{\! q}) & <
&8(\log_2q+2)+9,\\
\frac{1}{\log_23}\log_2q & \leqslant & T(Q_{4n}\times
\mathbb{Z}/_{\! q}) & < & 6\log_2q+18.
\end{array}$$
\item Let $n,h,s$ be integers, with $h,n\geqslant 3$ and $n,s$
coprime, and $s$ odd. Let $q=2^{h-2}s$ and suppose that $(n,q)$ is
a Zaremba pair. Then:
$$\begin{array}{rcccl}
\log_2s+h & < & c(D_{2^hn}\times\mathbb{Z}/_{\! s}) & < &
8(\log_2s+h)+15,\\
\frac1{\log_23}\log_2s & \leqslant &
T(D_{2^hn}\times\mathbb{Z}/_{\! s}) & < & 6(\log_2s+h)+6.
\end{array}$$
\end{enumerate}
In both the estimates on the complexity of the group, the lower
bound equals the base-$2$ logarithm of the order of the torsion of
the Abelianization, and the same term appears in the upper bound
too.
\end{theorem}

\begin{proof}
To begin, we notice that if $(n,q)$ is a Zaremba pair then $n<6q$.
Combining this inequality with the estimates obtained in
Theorem~\ref{Milnor:groups:thm}, we get the upper bounds on
complexity. The lower bounds on the complexity and the
$T$-invariant come from Theorems~\ref{lower:bound:thm}
and~\ref{Delzant:estimate:2:thm} respectively.

To get the upper estimates on the $T$-invariant, we first prove
the following assertion, first remarked by Delzant himself: {\em
If $M$ is a closed $3$-manifold then $T(\pi_1(M))$ does not exceed
twice the number of $3$-simplices in any triangulation of $M$}. To
see this, we note that a triangulation gives in particular a cell
decomposition of $M$, so we can employ the general algorithm
yielding a presentation of the fundamental group. The relations
then correspond to faces, and each face is a triangle, therefore
the length of any relation is no more than 3. Since $M$ is a
closed $3$-manifold, the number of faces is twice the number of
3-simplices, and the conclusion follows.

The fact just stated also holds for ``singular triangulations,''
{\em i.e.} triangulations with multiple and self-adjacencies.
Therefore it follows
from~\cite[Theorem 2.2.4]{Ma03} that if $M$ is irreducible and
$c(M)>0$ then $T(\pi_1(M))\leqslant 2c(M)$. Now we note that the
groups in the statement occur as fundamental groups of Seifert
manifolds of type $(S^2;(2,1),(2,1),(n,q),-1)$ (see, for
instance,~\cite[Chapter 2]{Ma03}). In particular, for odd $q$ this
manifold has fundamental group $Q_{4n}\times \mathbb{Z}/_{\! q}$,
and for even $q$ it has fundamental group $D_{2^hn}\times
\mathbb{Z}/_{\! {s}}$ with $h,s$ as in the statement. The
conclusion then easily follows from
Theorem~\ref{compl:mfds:thm}.\end{proof}

Using Theorem~\ref{complexity:Seifert:thm}, an argument similar to
that used in the proofs of Propositions~\ref{Petronio:prop}
and~\ref{Petronio:Milnor:prop} can now be employed to establish
the following:

\begin{proposition}\label{Petronio:Seifert:prop}
Let $G$ be one of the groups of
Theorem~\ref{complexity:Seifert:thm}, and let $\langle a_1,\ldots,
a_n|\ r_1,\ldots, r_m \rangle$ be a length-minimizing presentation
of $G$. Then
$$\frac{|r_1|+\ldots+|r_m|}{m}< 64.$$
\end{proposition}

\vspace{1.5cm}

\noindent
Chelyabinsk State University\\
ul. Br. Kashirinykh, 129,\\
454021 Chelyabinsk, Russia\\
pervova@csu.ru

\vspace{.5cm}

\noindent
Dipartimento di Matematica Applicata\\
Universit\`a di Pisa\\
Via Bonanno Pisano 25B, 56126 Pisa, Italy\\
petronio@dm.unipi.it

\end{document}